\documentclass[12pt]{article}
\usepackage{amssymb}
\usepackage{amsmath}

\usepackage{mathrsfs}

\usepackage{geometry}          

\usepackage{epstopdf}
\usepackage{amsthm}
\usepackage{enumerate}

\usepackage[]{tikz,pgfplots}
\usetikzlibrary{calc}

\newtheorem{theorem}{Theorem}

\newtheorem{proposition}[theorem]{Proposition}

\newtheorem{remark}{Remark}

\newcommand{\sgn}{\mbox{\rm sgn}}							% SIGN
\newcommand{\PP}{\mathbb{P}}                                                          % michal's standard shortcuts 
\newcommand{\RR}{\mathbb{R}}                                                           %  double letter instead of "mathbb"
\newcommand{\NN}{\mathbb{N}}                                                           
                                                           
\newcommand{\bc}{\mathcal B}                                                               % letter followed by "c" means the caligraphic version 

\def\qed{\unskip\nobreak\hfill\penalty50\hskip 3pt\hbox{}\nobreak
\hfill\hbox{\vrule width 4 pt height 10 pt}}

\title{Example of a non-standard Extreme Value Law}
\date{\today}
\begin{document}
\maketitle
\author{N Haydn 
\footnote{\tiny Department of Mathematics, University of Southern California,
Los Angeles, 90089-2532. E-mail: {\tt nhaydn@usc.edu}}
%\thanks{Mathematics Department, USC,
%Los Angeles, 90089-2532. E-mail: $<$nhaydn@usc.edu$>$.}
, M Kupsa
\footnote{\tiny Institute of Information Theory and
  Automation, The Academy of Sciences of the Czech Republic, Prague 8, CZ-18208. E-mail: {\tt kupsa@utia.cas.cz}}
\footnote{\tiny Faculty of Information Technology,
  Czech Technical University in Prague, Prague 6, CZ-16000.}
}
% \begin{center}{\tiny Department of Mathematics, University of Southern California,
% Los Angeles, 90089-2532. E-mail: {\tt nhaydn@usc.edu}.}\end{center}

\begin{abstract}
It has been shown that sufficiently well mixing dynamical systems with positive 
entropy have extreme value laws which in the limit converge to one of the 
three standard distributions known for i.i.d. processes, namely
Gumbel, Fr\'echet and Weibull distributions. In this short note we give an example
which has a non-standard limiting distribution for its extreme values.
Rotations of the circle by  irrational numbers are used and it will
be shown that the limiting distribution is a step function where the limit
has to be taken along a suitable sequence given by the convergents.
\end{abstract}

\section{Introduction}

For i.i.d.\ random variables extreme value theory is a well established topic. 
Given a sequence of random variables $X_j$ one forms the maximal random 
variable $M_n=\max_{1\le j\le n}X_j$. If there exists a sequence of numbers, $a_n, b_n$ 
such that $a_n(M_n-b_n)$ converges in distribution to a limit, then one says that
the sequence $(X_j)_{j\in\NN}$ has the extreme value property. It is know that 
for i.i.d. random variables $X_j$ satisfying some weak condition the 
limiting distribution is one of the three: Gumbel (the distribution function
is $e^{-e^{-t}}$), Fr\'echet (with the distribution function $e^{-t^{-\beta}}$ for $t>0$ and a parameter $\beta>0$)
or Weibull (with the distribution function $e^{-(-t)^\beta}$ for $t<0$ and a parameter $\beta>0$).

In a dynamical system with map $T$ on a space $\Omega$ the random variable chosen is a given 
function $X_0$ evaluated along the orbit, that is $X_j=X_0\circ T^j$. The function
chosen is rotational symmetric: a base point $x$ is selected and then $X_0(y)=g(d(x,y))$,
where $g$ is a function defined on $\mathbb{R}^+$. Typically $g(s)=-\log s$ is chosen.
The pursuit of extreme values in dynamics is quite recent with the most notable first 
paper in this domain being~\cite{Col} where for non-uniformly expanding $C^2$ maps on 
the interval the limiting extreme value law for $g(s)=-\log s$ was proven to be Gumbel.
For the quadratic map on the interval and `Benedicks-Carleson' parameters it was shown in~\cite{FF2}
that the extreme value statistics for $g(s)=C-s^\beta$, for a constant $C>0$ and a parameter $\beta$,
tend in the limit to a Weibull distribution.
For more general non-uniformly hyperbolic maps the EVL statistics was addressed in~\cite{HNT}. 
For higher dimensional hyperbolic maps with discontinuities, like dispersing billiards, Lozi map
and Lorenz type maps the limiting EVL was established in~\cite{GHN}.
Since the functions $g$ are always connected to the metric, the fact that $M_n$ is large 
for some time geometrically means that a point does not enter the neighborhood of $x$
for this time. The EVL property and the distribution of hitting or return times are therefore
intimately connected to each other. In fact the equivalence was formally established 
first for absolutely continuous measures in~\cite{FFT1} and then for more general measures
in~\cite{FFT2}.

Here we provide an example that is contrary to all the quoted results. It yields a non-standard
limiting EVL. Since we use circle rotations, we don't obtain the good mixing properties of the 
systems mentioned nor a good decay of correlations. In fact, it is known~\cite{KP} that
the hitting times distribution for circle maps are `non-standard'. The limiting distributions
turn out to be locally constant and not exponential as is witnessed by many hyperbolic
systems or those that at least display a sufficiently fast decay of correlations
(see e.g.~\cite{Aba1, Hay}).

\section{General settings}
Let us consider a ``probability'' dynamical system $(\Omega,\bc,\PP,T)$
with invariant measure $\PP$ and a measurable function $X$ from $\Omega$ to $\RR$. Denote the maximum value in the first $n$ trials by $M_n$, i.e. $M_n=\max_{1\le j\le n}X_j$, where $X_j(x)=X\circ T^j$. The question is whether there exists sequences $(a_n)$ and $(b_n)$ such that the rescaled random variables $a_n(M_n-b_n)$ converges in distribution and what the limit is. 

Let us recall that real random variables $Y_n$, $n\in\NN$, converge to a real random variable $Y$ in distribution if and only if the distribution functions $F_{X_n}(t)=\PP(X_n\leq t)$ converges to the distribution function $F_X(t)=\PP(X \leq t)$ at every point $t$, where the limiting distribution $F_X$ is continuous. It is a simple observation, that the convergence in distribution can be expressed in the same way in the terms of functions $\PP(X_n> t)$ and $\PP(X>t)$. 

Hence, we ask if the functions
$$H_n(y):=\PP\left((M_n-b_n)a_n>y\right), y\in\RR$$
converge to a right-continuous decreasing (not necessarily strictly) function $H(y)$ at every point $y\in\RR$, where the function $H$ is continuous. Since any right-continuous function has at most countably many points of discontinuity, the limiting distribution $H$ is uniquely determined (if exists). It represents a real random variable if the limits of $H(y)$ at plus and minus infinity are $1$ and $0$, respectively.

The maximum value statistics are tightly connected with the statistics of the entry times. For a set $B\in \bc$ denote by $\tau_B$ the entry time function which is 
 given by 
 $$
 \tau_B(x)=\min\{j\ge1: T^j(x)\in B\}
 $$
 ($\tau_B(x)=\infty$ if $x$ never enters $B$).
 The normalized entry times distribution then is given by
 $$
 F_B(t)=\mathbb{P}\left(\tau_B\le \frac{t}{\PP(B)}\right)
 $$
 for $t\in\RR$, is increasing, constant 
 on intervals $[j\PP(B),(j+1)\PP(B))$ for $j=0,1,\dots$ and the jumps at the points 
 $j\PP(B)$, $j=1,2,\dots$, are less or equal to $\PP(B)$.

The connection between the maximum value $M_n$ and the entry times can be expressed in terms of level sets:
\[   L(y)=\{x\in\Omega, X(x)>y\},\qquad y\in\RR. \]
For the distribution of the maximal values we have 
$$
\PP(M_n>y)=\PP(\tau_{L(y)}\leq n),
$$
since the sets on the both sides equal.
% Let us make the following assumption on the existence of a limit law for the entry times:
% \begin{itemize}
% \item For $y\in\RR$ let 
% $$
%  F_{L(y)}(t)=\PP\left(\PP(L(y))\tau_{L(y)}\leq t\right),\qquad t>0
%  $$
%  be the entry distribution function of $L(y)$  and assume that it converges in distribution to a limit law 
%  $F$, when $y$ is going to infinity, i.e. $\lim_{y\to\infty}F_{L(y)}(t)=F(t)$ for every $t\in\RR$.
% \end{itemize}

Consequently, for sequences $a_n, b_n, n=1,2,\dots$, the rescaled maximum value variables $H_n$ satisfy the following equalities:
\begin{align}\label{eq:general-max-ret}
H_n(y)&:=\PP\left((M_n-b_n)a_n>y\right)=\PP\left(M_n>\frac1{a_n}y+b_n\right)=\PP\left(\tau_{L(\frac{y}{a_n}+b_n)}\le n\right)\\
&=F_{L\left(\frac{y}{a_n}+b_n\right)}\left(n\ \PP\left(L\left(\frac{y}{a_n}+b_n\right)\right)\right).
\end{align}
We will use this equality in the next section where we calculate directly the limiting distribution for the sequence $H_n$ in irrational rotation of the interval. 
% In Section \ref{sec:rot2}, we derive the limiting form of the equality above and use it to calculate again the same limiting distribution using the limit entry times law.

% To get a limit law for the rescaled maximum value it is enough if the sequences $(a_n)$ and $(b_n)$ satisfy the following conditions:
% \begin{itemize}
% \item $\frac{y}{a_n}+b_n$ converges to infinity for every $y\in \RR$,
% \item $n\cdot \PP\left(L\left(\frac{y}{a_n}+b_n\right)\right)$ converges for every $y\in \RR$ 
%  to a function $g(y)$.
% \end{itemize}
% Then
% \[H(y):=\lim_{n\to\infty}H_n(y)=F(g(y)).\]   

\section{Rotation of  the interval}

Let us consider a rotation $T:[0,1)\rightarrow [0,1)$,  $Tx=x+\alpha \mbox{ mod } 1$, on the unit-interval (or circle) by an irrational angle $\alpha\in(0,1)$. The Lebesgue measure $\mu$ is then the only invariant probability measure. We consider the continued fraction expansion
$$
\alpha=[c_1,c_2,c_3,\dots]=\cfrac{1}{c_1+\cfrac1{c_2+\cfrac1{c_3+\ldots}}}.
$$
The convergents of $\alpha$ are then $\frac{p_k}{q_k}$ where $p_k=c_kp_{k-1}+p_{k-2}$, $p_0=0,p_1=1$
 and  $q_k=c_kq_{k-1}+q_{k-2}$, $q_0=1,q_1=c_1$. The numbers $q_k\alpha-p_k$, $k\in\NN$, 
 form an alternating positive and negative sequence and their absolute values $\eta_k=|q_k\alpha-p_k|$ satisfy the implicit formula $\eta_k=\eta_{k-2}-c_k\eta_{k-1}$, $\eta_0=\alpha$, $\eta_1=1-c_1\alpha$.

Denote the following nested sequence of intervals, 
$$
 B_k=
 \begin{cases}
   (-\eta_{k+1},\eta_k)&\mbox{for $k$ even}\\
   (-\eta_k,\eta_{k+1}) &\mbox{for $k$ odd}.
 \end{cases}
$$
With no danger of ambiguity, the sets $B_k$, $k\in\NN$, are considered as subsets of the state space $[0,1)$, so we identify the above-mentioned intervals  with their images under the projection $\text{ mod }1:\mathbb{R}\to [0,1)$.

Let $X$ be a random variable on $[0,1)$ defined as follows:
$$
%X(x)=q_{k+1},\quad x\in B_k\setminus B_{k+1},\qquad X(x)=q_0\quad x\not\in B_0.
X(x)=q_{k+1},\qquad \text{where } k=\min\{\ell\in\NN\mid x\not\in B_\ell\}.
$$ 

%%%%%%%%%%%%%%%%%%%%%%%%%%%%%%%%%%%%%%%%%%%%%%%
%%%%%%%    FIGURE 1
\begin{figure}[htbp]
  \centering

\begin{tikzpicture}[yscale=0.5, xscale=14]
%  \draw[very thin,color=gray] (-3.1,-4.1) grid (3.9,3.9);
  \pgfmathsetmacro{\gold}{1.618}; % golden number  
  \pgfmathsetmacro{\invgold}{0.618}; % reciprocal for golden number
  \pgfmathsetmacro{\betaa}{0.723}; % (a+a^2)/(1+a^2), where a is \invgold
  \pgfmathsetmacro{\betab}{1.17}; %  (a+a^2)/(1+a^2)x1/a, where a is \invgold

\pgfmathsetmacro{\invggold}{\invgold*\invgold};
\pgfmathsetmacro{\invgggold}{\invggold*\invgold};
\pgfmathsetmacro{\invggggold}{\invgggold*\invgold};
\pgfmathsetmacro{\invgggggold}{\invggggold*\invgold};
\pgfmathsetmacro{\invggggggold}{\invgggggold*\invgold};
\pgfmathsetmacro{\invgggggggold}{\invggggggold*\invgold};
\pgfmathsetmacro{\invggggggggold}{\invgggggggold*\invgold};

  \draw[->] ($(\invgold-1.05,0)$)  -- ($(\invgold+0.05,0)$) node[below] {$x$};
  \draw[->] (0,-0.2) node[below] {$0$} -- (0,8) node[left] {$X(x)$};
  
\draw[-] (\invgold,-0.2) node[below] {$\eta_{2k}$} -- (\invgold,0.2);
\draw[-] ($(-\invggold,-0.2)$) node[below] {$-\eta_{2k+1}$} -- ($(-\invggold,0.2)$);
\draw[-] (\invgggold,-0.2) node[below] {$\eta_{2k+2}$} -- (\invgggold,0.2);
\draw[-] ($(-\invggggold,-0.2)$) node[below] {$-\eta_{2k+3}$} -- ($(-\invggggold,0.2)$);
\draw[-] (\invgggggold,-0.2) node[below] {$\eta_{2k+4}$} -- (\invgggggold,0.2);
\draw[-] ($(-\invggggggold,-0.2)$) node[below] {$-\eta_{2k+5}$} -- ($(-\invggggggold,0.2)$);

\draw[-,thick] ($(\invgggold,1)$)  -- node[above] {$q_{2k+1}$} ($(\invgold,1)$);
\draw[-,thick] ($(-\invggold,\gold)$)  -- node[above] {$q_{2k+2}$} ($(-\invggggold,\gold)$);
\draw[-,thick] ($(\invgggggold,1/\invggold)$)  -- node[above] {$q_{2k+3}$} ($(\invgggold,1/\invggold)$);
\draw[-,thick] ($(-\invggggold,1/\invgggold)$)  -- node[above] {$q_{2k+4}$} ($(-\invggggggold,1/\invgggold)$);
\draw[-,thick] ($(\invgggggggold,1/\invggggold)$)  -- node[above] {$q_{2k+5}$} ($(\invgggggold,1/\invggggold)$);
% \draw[-,thick] (\ggold,\betaa/\ggold)  -- (\gggold,\betaa/\ggold);
% \draw[-,thick] (\gggold,\betaa/\gggold)  -- (\ggggold,\betaa/\gggold);
% \draw[-,thick] (\ggggold,\betaa/\ggggold)  -- (\rend,\betaa/\ggggold);

\end{tikzpicture}

  \caption{Initial distribution $X$}
  \label{fig:inidist}
\end{figure}
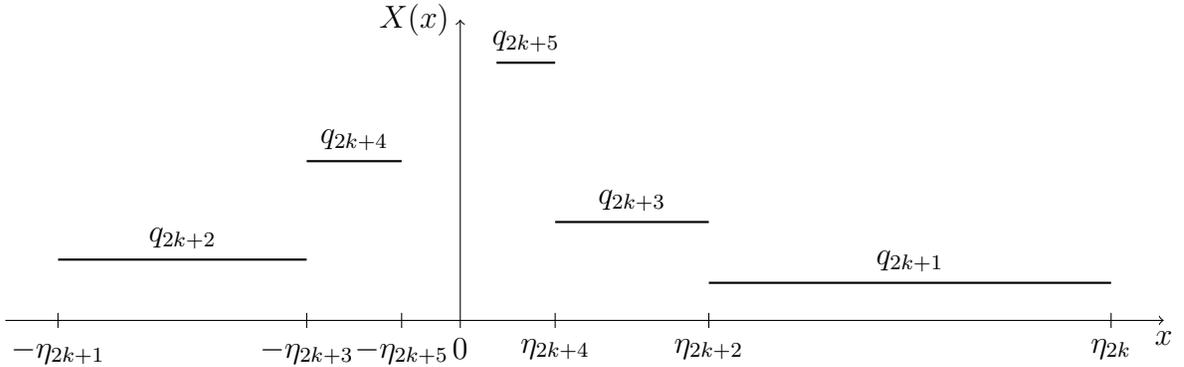

The sets $B_k$ were chosen to have nice return times, namely:
\begin{align*}
  \tau_{B_k}(x)=
  \begin{cases}
    q_{k+1}&\mbox{if $x\cdot\sgn(q_k\alpha-p_k)\in (0,\eta_k)$}\\
    q_{k} &\mbox{if $x\cdot\sgn(q_k\alpha-p_k)\in (-\eta_{k+1},0)$}.
  \end{cases}
\end{align*}

Consequently we obtain that the distribution of the entry time satisfy the following conditions (for the detailed proof see Proposition 3 in \cite{CF}): 
\begin{align}\label{eq:entry times to B_k}
  \mu\left(\tau_{B_k}\le s\right)=
  \begin{cases}
    1&\mbox{if $s\ge q_{k+1}$}\\
    q_{k}\eta_{k+1}+\eta_k[s]&\mbox{if $s\in[q_{k},q_{k+1})$}\\
    ( \eta_k+\eta_{k+1})[s]&\mbox{if $s\in[0,q_{k})$.}
  \end{cases}
\end{align}
   (note that $q_{k+1}\eta_k+\eta_{k+1}q_{k}=1$), since
   $q_k(\eta_k+\eta_{k+1})+\eta_k[s-q_k]=q_{k}\eta_{k+1}+\eta_{k+1}[s]$.

For an infinite set of integers $K\subset\NN$, denote the following limits (if they exist):
\begin{align*}
  \lim_{k\in K}\frac{q_{k+j}}{q_k}=\gamma_j,\quad\lim_{k\in K}\frac{\eta_{k+j}}{\eta_k}=\delta_j, \qquad \forall j\in\NN.
\end{align*}

It is a standard fact that $q_{k+2}/q_k>2$ for every $k$. Hence, by the definition, the sequence $\gamma_j$ is increasing (not necessarily strictly) and goes to infinity, whenever the limits $\gamma_j$'s exist. Moreover, 
\begin{align}\label{eq:qeta}
\lim_{k\in K} q_k\eta_k=\frac{1}{\gamma_1+\delta_1},
\end{align}
whenever $\gamma_1$ and $\delta_1$ exist. Indeed,
\begin{align*}
\lim_{k\in K} q_k\eta_k&=\lim_{k\in K} \left(\frac{1}{q_k\eta_k}\right)^{-1}=\lim_{k\in K} \left(\frac{q_{k+1}\eta_k+q_k\eta_{k+1}}{q_k\eta_k}\right)^{-1}\\
&=\lim_{k\in K} \left(\frac{q_{k+1}}{q_k}+\frac{\eta_{k+1}}{\eta_k}\right)^{-1}=(\gamma_1+\delta_1)^{-1}.
\end{align*}

In terms of these limits, we express our main theorem.

%%%%%%%%%%%%%%%%%%%%%%%%%%%%%%%%%%%%%%%%%%%
%%%%%%%%%%%%%%%%%        MAIN THEOREM
\begin{theorem}\label{thm:rotations} Let $K$ be an infinite set of natural numbers and the limits $\gamma_j$ and $\delta_j$ exist (along $K$) for every $j\geq 0$. Then the random variables $M_{q_k}/q_k$, $k\in K$, converge in distribution to a random variable $M$ with the following distribution: 
$$
H(y)=\PP(M>y)=
\begin{cases}
  1&\mbox{if $y<1$}\\
  \frac{\delta_j+\delta_{j+1}}{\gamma_1+\delta_1}&\mbox{if $y\in
    [\gamma_j,\gamma_{j+1})$, $ j=0,1,2,\dots$.}
\end{cases}
$$
\end{theorem}

\noindent{\bf Proof.}
Let us recall a standard fact, that convergence in distribution of the random variables $M_{q_k}/q_k$ to a random variable $M$ agrees with the pointwise convergence of distribution functions $H_{q_k}$ to $H$ at every point except the points of discontinuities of the limiting function $H$. By the formula for $H$, its points of discontinuities are the points $\gamma_j$, $j\in\NN$. In addition, $\gamma_j$ goes to infinity (or reaches infinity in a finite step). Thus, we need to prove that $H_{q_k}(y)$ tends to $H(y)$ for every $y$ from the intervals $(-\infty,\gamma_1)$ and $(\gamma_j,\gamma_{j+1})$, $j\in\NN_0$. 
We treat the two cases separately as follows:\\
{\bf (I)} Let $y\leq \gamma_0=1$. Then for every $k\in K$, $yq_k< q_k$, $B_{k-1}\subset L(yq_k)$ and
$$
H_{q_k}(y)=\mu\left(M_{q_k}/q_k>y\right)=\mu\left(\tau_{L(q_ky)}\le q_k\right)\geq \mu\left(\tau_{B_{k-1}}\le q_k\right)=1.
  $$
The last equality follows from (\ref{eq:entry times to B_k}). Hence the limit $\lim_{k\in K} H_{q_k}(y)$ is one.\\
{\bf (II)} Let $y\in (\gamma_j,\gamma_{j+1})$ for some $j\geq 0$. By the definition of $\gamma_j$'s, $q_{k+j}< q_ky<q_{k+j+1}$ eventually for every $k\in\NN$. Hence $L(q_ky)=B_{k+j}$ for $k$ big enough. It implies that
$$
H_{q_k}(y)=\mu\left(\tau_{B_{k+j}}\le q_k\right)=q_k(\eta_{k+j}+\eta_{k+j+1}).
$$
The last equality follows from (\ref{eq:entry times to B_k}) and the fact that $q_k<q_{k+j}$. By \eqref{eq:qeta},
$$
H_{q_k}(y)=q_k\eta_k\left(\frac{\eta_{k+j}}{\eta_k}+\frac{\eta_{k+j+1}}{\eta_k}\right) \rightarrow\frac{\delta_j+\delta_{j+1}}{\gamma_1+\delta_1}
.$$
\qed

\begin{remark}
As we mentioned above Theorem \ref{thm:rotations}, the sequence $(\gamma_j)_{j\in\NN}$ is increasing (not necessarily strictly) and goes to infinity. Hence, the description of $H(y)$ from Theorem \ref{thm:rotations} characterizes the distribution function in a correct and unique way. Let us notice, that $(\gamma_j)_{j\in\NN}$ can be eventually equal to $+\infty$. We discuss this case in Section \ref{sec:finite-lim-gamma}.
\end{remark}

\section{The limits $\gamma_j$ and $\delta_j$}
\label{sec:discuss existence of limits}
A natural question is, for a given irrational rotation, what is the set $K$ and what are the limits $\gamma_j$ and $\delta_j$. Let the angle $\alpha$ be given in the form of continued fraction $\alpha=[c_1,c_2,\ldots]$, 
and assume $K\subset\NN$ is infinite. 
%% The number $[c_1,c_2,\ldots]$ is represented by the point $(c_n)$ in the product of discrete topological spaces $\NN^\NN$. Before we discuss particular cases, let us remark, that if the limits $\gamma_j$ and $\gamma_j'$ exists, then 
Put
\[
\nu_j=\lim_{k\in K}\frac{q_{k+j-1}}{q_{k+j}},\qquad 
\theta_j=\lim_{k\in K}\frac{\eta_{k+j}}{\eta_{k+j-1}},\qquad j\geq 1,
\]
It is a standard fact that
\[
\frac{q_{k-1}}{q_k}=[c_k,c_{k-1},\ldots,c_1],\qquad\frac{\eta_k}{\eta_{k-1}}=[c_{k+1},c_{k+2},\ldots],\qquad k\in\NN.
\]
Hence, the limits $\nu_j$ and $\theta_j$ can be expressed in another way:
\[
\nu_j=\lim_{k\in K} [c_{k+j},c_{k+j-1},\ldots,c_1],\qquad \theta_j=\lim_{k\in K} [c_{k+j+1},c_{k+j+2},\ldots],\qquad j\geq 0,
\]

It follows from the definition that for all $j\ge1$:
\[
\gamma_j=\frac{\gamma_{j-1}}{\nu_j}=\prod^j_{i=1}(\nu_i)^{-1},\qquad \delta_j=\theta_j\delta_{j-1}=\prod^j_{i=1}\theta_i,
\]
whenever the limits used in the equalities exist (where $0^{-1}=\infty$).

\subsection{Finite limits $\gamma_j$}
\label{sec:finite-lim-gamma}
First, let us suppose that the limits $\gamma_j=\lim_{k\in K}\frac{q_{k+j}}{q_k}$ exist and are finite for every $j\geq 1$. This condition ensures that
\begin{itemize}
\item the limits $\nu_j$, $j\geq 1$, exist and are nonzero,
\item the limits  $\lim_{k\in K} c_{k+j}$, $j\geq 1$, exist and are finite, i.e. for every $j\geq 1$, the sequence $(c_{k+j})_{k\in K}$ is eventually constant, %%% this is commented because it is not equivalent, it is weaker condition.
\item the limits $\theta_j$, $j\geq 1$, exist and are nonzero,
\item the limits  $\delta_j$, $j\geq 1$, exist and are nonzero.
\end{itemize}
In such a settings the limit distribution function for extremes $H(y)$ described in the main theorem has {\bf countably} many jumps (down) in the points $\gamma_j$, $j\geq 0$. Letting $y$ go to$+\infty$, the function converges to $0$, but never reach this value.

\subsection{Infinite limits $\gamma_j$}
\label{sec:finite-lim-gamma}
% In this case we obtain that the limiting distribution $H$ is trivial, namely equal to the characteristic
% functions  $1_{[-\infty,1)}$ and the same is true for the survival function of the constant random 
% variable $M=1$ .  
Let us suppose that the limits $\gamma_j=\lim_{k\in K}\frac{q_{k+j}}{q_k}$ exist for every $j\geq 1$ and some of them are infinite. In this case, the situation is more complex. We look at the two cases when $\gamma_1$ is finite and infinite separately as follows:\\
{\bf (I) $\gamma_1<\infty$:}
Let $N$ be such an index, that $\gamma_N$ is the first infinite member of the sequence $(\gamma_j)_{j\geq 1}$. Then $N\geq 2$ and the following conditions hold:
\begin{itemize}
\item the limits $\nu_j$, $1\leq j< N$, exist and are nonzero. The limit  $\nu_N$ exists and is zero.
\item The limits  $\lim_{k\in K} c_{k+j}$, $1\leq j< N$, exist and are finite, i.e. for every $1\leq j< N$, the sequence $(c_{k+j})_{k\in K}$ is eventually constant. The limit $\lim_{k\in K} c_{k+N}$ is infinite. %%% this is commented because it is not equivalent, it is weaker condition.
\item The limits $\theta_j$, $1\leq j< N-1$, exist and are nonzero. The limit  $\nu_{N-1}$ exists and is zero.
\item The limits  $\delta_j$, $1\leq j< N-1$, exist and are nonzero. The limit $\delta_{N-1}$ exists and is zero.
\end{itemize}
In such a setting the limit distribution function for extremes $H(y)$ described in the main theorem has {\bf finitely} many jumps (down) at the points $\gamma_j$, $0\leq j\leq N-1$. The function reaches zero at the point $\gamma_{N-1}$, indeed, 
\[
H(y)=\frac{\delta_{N-1}+\delta_N}{\gamma_1+\delta_1}=0,\qquad\text{ for } y\in[\gamma_{N-1},\gamma_N)=[\gamma_{N-1},\infty). 
\]

\noindent {\bf (II) $\gamma_1=\infty$:}\label{rem:trivial}
In this case, the limiting distribution $H$ is the characteristic function  $1_{[-\infty,1)}$ which
 can easily be seen following the steps in the proof of Theorem \ref{thm:rotations}, even for the case when the limits $\delta_j$, $j\geq 1$, do not exist.  
 
 \vspace{3mm}

 \noindent Theorem \ref{thm:rotations} ensures $H$ to be non-trivial in all cases, except two: when $\gamma_1$ is infinite, or when $\gamma_1=1$ and $\gamma_2$ is infinite. The former case was already discussed above. In the latter case, the points of discontinuities, $\gamma_0$ and $\gamma_1$, are equal. Thus, there is only one point of the discontinuity at $\gamma_0=1$. As we already mentioned above, if $\gamma_2$ is infinite, then $\delta_1$ is zero, so is $\delta_2$. We get, that the value of the function $H$ on the interval $[1,\infty)=[\gamma_1,\gamma_2)$ is zero.

\vspace{3mm}

\noindent {\bf Example.} The two cases when the limiting law for the extreme values is trivial, as well as the case of the non-trivial case with a finite number of jumps is illustrated by the rotation numbers
$$
\alpha=[\overbrace{1,1,\ldots,1}^{N-1},2,\overbrace{1,1,\ldots,1}^{N-1},3,\overbrace{1,1,\ldots,1}^{N-1},4,
\dots].
$$
 where $N\ge1$. Here we choose $K=\{kN: k\in\mathbb{N}\}$. 

If $N\ge2$, we obtain that 
 $$
 \frac{q_{kN+j-1}}{q_{kN+j}}=[\overbrace{1,\ldots,1}^{j},k+1,\overbrace{1,\ldots,1}^{N-1},k,
 \overbrace{1,\ldots,1}^{N-1},k-1,\dots]
$$
for every $j=0,\dots,N-1$. In particular we immediately get $\nu_0=0$, since
$$
\nu_0=\lim_{k\in\NN} [k+1,\overbrace{1,\ldots,1}^{N-1},k,
\overbrace{1,\ldots,1}^{N-1},k-1,\dots]\leq\lim \frac{1}{k+1}=0.
$$
For $j=1,\dots,N-1$ we obtain
\begin{align*}
\nu_{j}&=\lim_{k\in\NN} [\overbrace{1,\ldots,1}^{j},k+1,\overbrace{1,\ldots,1}^{N-1},k,
\overbrace{1,\ldots,1}^{N-1},k-1,\dots]\\
&=\lim_{k\in\NN}
\frac{1}{1+[\overbrace{1,\ldots,1}^{j-1},k+1,\overbrace{1,\ldots,1}^{N-1},k,
  \overbrace{1,\ldots,1}^{N-1},k-1,\dots]}=\frac{1}{1+\nu_{j-1}},
\end{align*}
This recursive formula ensures that $\nu_j=\frac{s_{j}}{s_{j+1}}$, for every $j=0,\ldots,N-1$, where $s_j$ is the Fibonacci sequence given recursively by
 $s_{j+1}=s_j+s_{j-1}$ and $s_0=0$, $s_1=1$. For $j=N$ we obtain $\nu_{N}=\nu_0=0$. 
  Hence $\gamma_j=s_{j+1}$ for $j=1,\ldots,N-1$ and $\gamma_{N}=\infty$.
  On the other hand
  $$
 \frac{\eta_{kN+j}}{\eta_{kN+j-1}}=[\overbrace{1,\ldots,1}^{N-1-j},k+2,\overbrace{1,\ldots,1}^{N-1},k+3,
 \overbrace{1,\ldots,1}^{N-1},k+4,\dots],
  $$
for every $j=0,\ldots,N-1$. For the value $j=N-1$ we immediately get  $\theta_{N-1}=\nu_0=0$.
For $j=2,\ldots,N-1$ one obtains $\theta_j=\nu_{N-1-j)}=\frac{s_{N-1-j}}{s_{N-j}}$ and
therefore $\delta_j=s_{N-j-1}/s_{N-1}$, for every $j=1,\ldots,N-1$. In particular, $\delta_{N-1}=0$. It implies that function $H$ vanishes on $[s_N,\infty)$.

If $N\geq 3$, the function $H$ is non-trivial and has $N-1$ jumps at the points $\gamma_j=s_{j+1}$, $j=1,\ldots, N-1$. By Theorem \ref{thm:rotations}, for every $j=1,\ldots,N-2$, the value of $H$ on the interval $[s_{j+1},s_{j+2})$ is equal to $s_{N-j}/s_{N}$. The function is zero on $[s_N,\infty)$. 

If $N=2$, then $\gamma_1=s_2=1$ and $\gamma_2$ is infinite. Here we get the trivial limit distribution $H$.

  If $N=1$ then
  $$
\alpha=[2,3,4,\dots].
$$
and we get $\nu_1=0$ which implies $\gamma_1=\infty$. Here  $K=\mathbb{N}$.
This is the second case when  $H$ is the trivial limiting distribution.

If the sequence of entries $2,3,4,\dots$ in the continued fraction expansion are replaced by a 
sequence of numbers $n_1,n_2, n_3,\dots$ which converges to infinity then we will get the 
same values for  $\gamma_j$ and $\delta_j$. 

\vspace{3mm}
  
Below we apply Theorem~\ref{thm:rotations} to two classical situations: rotation numbers of 
constant type and diverging rotation numbers.

\subsection{Rotation numbers of constant type}
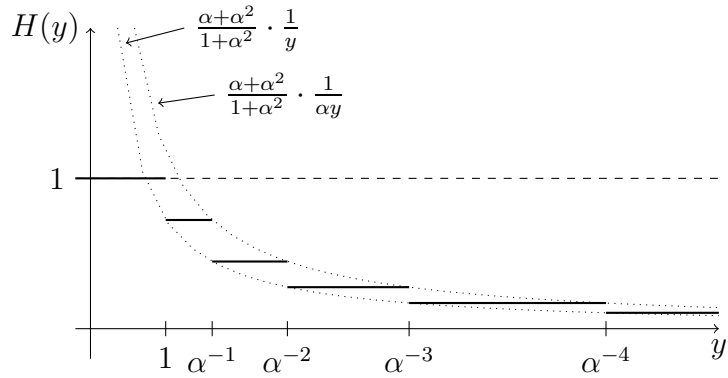
\begin{figure}[htbp]
  \centering

\begin{tikzpicture}[yscale=2]
%  \draw[very thin,color=gray] (-3.1,-4.1) grid (3.9,3.9);
  \pgfmathsetmacro{\gold}{1.618}; % golden number  
  \pgfmathsetmacro{\invgold}{0.618}; % reciprocal for golden number
  \pgfmathsetmacro{\betaa}{0.723}; % (a+a^2)/(1+a^2), where a is \invgold
  \pgfmathsetmacro{\betab}{1.17}; %  (a+a^2)/(1+a^2)x1/a, where a is \invgold

\pgfmathsetmacro{\ggold}{\gold*\gold};
\pgfmathsetmacro{\gggold}{\ggold*\gold};
\pgfmathsetmacro{\ggggold}{\gggold*\gold};

  \pgfmathsetmacro{\rend}{\ggggold+1.5}; % sirka
 
  \draw[->] (-0.2,0) -- (\rend,0) node[below] {$y$};
  \draw[-, dashed] (-0.2,1) node[left] {$1$} -- (\rend,1);
  \draw[->] (0,-0.2) -- (0,2) node[left] {$H(y)$};
  
  \draw[dotted,domain=0.5*\betaa:\rend] plot (\x,{\betaa/\x});
\node (firsthyp) at (\rend/4,2) {$\frac{\alpha+\alpha^2}{1+\alpha^2}\cdot \frac1y$};
\draw[->] (firsthyp.west) -- ($(0.65*\betaa, 1.9)$);
\node[below,anchor=north] (secondhyp) at ($(firsthyp.south)+(0.5,0)$) {$\frac{\alpha+\alpha^2}{1+\alpha^2}\cdot \frac1{\alpha y}$};
  \draw[dotted,domain=0.5*\betab:\rend] plot (\x,{1/\x*\betab});
\draw[->] (secondhyp.west) -- ($(0.9,1.5)$);

\draw[-] (1,-0.05) node[below] {$1$} -- (1,0.05);

\draw[-] (\gold,-0.05) node[below] {$\alpha^{-1}$} -- (\gold,0.05);
\draw[-] (\ggold,-0.05) node[below] {$\alpha^{-2}$} -- (\ggold,0.05);
\draw[-] (\gggold,-0.05) node[below] {$\alpha^{-3}$} -- (\gggold,0.05);
\draw[-] (\ggggold,-0.05) node[below] {$\alpha^{-4}$} -- (\ggggold,0.05);

\draw[-,thick] (-0.2,1)  -- (1,1);
\draw[-,thick] (1,\betaa)  -- (\gold,\betaa);
\draw[-,thick] (\gold,\betaa/\gold)  -- (\ggold,\betaa/\gold);
\draw[-,thick] (\ggold,\betaa/\ggold)  -- (\gggold,\betaa/\ggold);
\draw[-,thick] (\gggold,\betaa/\gggold)  -- (\ggggold,\betaa/\gggold);
\draw[-,thick] (\ggggold,\betaa/\ggggold)  -- (\rend,\betaa/\ggggold);

\end{tikzpicture}

  \caption{Limiting extreme value law for the rotation numbers of constant type.%, see Corollary \ref{cor:constant}.
}
  \label{fig:retlim}
\end{figure}

In this case one has $\alpha=[c,c,c,\ldots]$ that is $\alpha=\frac12(\sqrt{c^2+4}-c)$. The next corollary and Figure \ref{fig:retlim} describe the limiting extreme value law for these special numbers. 

In this case, the sequences $(\nu_j)_{j\in\NN}$ and $(\theta_j)_{j\in\NN}$ are constant, namely
\begin{align*}
\nu_j&=\lim_{k\in\NN}\frac{q_{k+j-1}}{q_{k+j}}=\lim_{k\in\NN}[\overbrace{c,c,\ldots,c}^{k+j}]=[c,c,\ldots]=\alpha,\\
\theta_j&=\lim_{k\in\NN} [c,c,\ldots]=[c,c,\ldots]=\alpha.
\end{align*}
Thus,
$$
\gamma_j=\alpha^{-j},\qquad \delta_j=\alpha^{j}, \qquad j\in\NN,
$$
Applying Theorem \ref{thm:rotations} we get that the random variables $M_{q_k}/q_k$, $k\in \NN$, converge in distribution to a random variable $M$ with the distribution: 
$$
H(y)=\PP(M>y)=
\begin{cases}
  1&\mbox{if $y<1$}\\
  \alpha^{j}\frac{\alpha+\alpha^2}{1+\alpha^2}&\mbox{if $y\in
    [\alpha^{-j},\alpha^{-j-1})$, $ j=0,1,2,\dots$.}
\end{cases}
$$

%%%%%%%%%%%%%%%%%%%%%%%%%%%%%%%%%%%%%%%%%%%%%%%
%%%%%%%    FIGURE 2

\vspace{3mm}

\noindent Let us note that in the case when $c=1$ we obtain the golden mean
 $\alpha=\frac{\sqrt5-1}2=[1,1,1,\dots]$ and for the limiting distribution
 (as $\alpha+\alpha^2=1$)
$$
H(y)=\PP(M>y)=
\begin{cases}
  1&\mbox{if $y<1$}\\
\frac{\alpha^j}{2-\alpha}&\mbox{if $y\in    [\alpha^{-j},\alpha^{-j-1})$, $ j=0,1,2,\dots$.}
\end{cases}
$$

\subsection{Divergent rotation number}
Let $\alpha=[c_1,c_2,c_3,\ldots]$ be such that $c_n$ converges to infinity. Let $K=\NN$.
In this case the coefficient $\gamma_1$ is equal to $+\infty$ and the random variables $M_{q_k}/q_k$ converge in distribution to the constant random variable $M=1$ (see Section~\ref{rem:trivial}). 

\section{Limiting behavior for rotations}
\label{sec:rot2}
In this section we show another way how to determine the limiting distribution of extreme values using the limiting distribution for the entry time.

Let the same parameters be as before $n=q_k$, $b_n=0$, $a_n=1/q_k$ and the infinite set $K\subset\NN$ according to the assumption of Theorem~\ref{thm:rotations}. 
Let us assume that the following two conditions hold for all $y\in\RR$:
\begin{itemize}
\item the following limit exists:
\begin{eqnarray*}
   g(y)&=&\lim_{k\in K} q_k \mu(L(q_ky))
\end{eqnarray*}
\item the measure of the level sets $\mu(L(q_ky))$ goes to zero, when $k$ goes to infinity,
\item the sequence of distribution functions for entry times $F_{L(q_ky)}$, $k\in K$, converges uniformly to a distribution function $\phi_y$.
\end{itemize}

Note that the last two conditions imply in particular that $\phi_y$ is continuous.

Under these three conditions, the extreme value distribution is
$$
H(y)=\lim_{k\in K} H_{q_k y}(y)=\lim_{k\in K} F_{L(q_ky)}(q_k \mu(L(q_ky))) = \phi_y(g(y)),
$$
where the last equality follows from the fact that $\phi_y$ need to be continuous (see \cite{KL}). The next proposition shows, that these assumptions are valid for every $y$ for which Theorem \ref{thm:rotations} ensures the existence of a nontrivial limit of $H_{q_k}(y)$.

%%%%%%%%%%%%%%%%%%%%%%%%%%%%%%%%%%%%%%%%%%%%%
%%%%%%%%%%%%     PROPOSITION
%%%%%%%%%%%%%%%%%%%%%%%%%%%%%%%%%%%%%%%%%%%%%
\begin{proposition}
Let $K$ be an infinite set of integers $K\subset\NN$, such that the limits $\gamma_j$ and $\delta_j$ exist for every $j\in\NN$. Assume  $\gamma_{j+1}<\infty$ for some $j\in\NN$.
Then for $y\in (\gamma_j,\gamma_{j+1})$ the following hold:\\
(i) The limit $\lim_{k\in K} \mu(L(q_k y))$ is zero,\\
(ii) The limit $g(y)=\lim_{k\in K} q_k\mu(L(q_ky))$ exists, is finite and satisfies the following equality:
$$g(y)=\frac{\delta_j+\delta_{j+1}}{\gamma_1+\delta_1}.$$
(iii) The distribution functions for entry times $F_{L(q_k y)}$, $k\in K$, converge uniformly to the distribution function $\phi_y$ that linearly interpolates the points
$$(0,0),\quad \left(\frac{(1+\theta_{j+1})\nu_{j+1}}{1+\theta_{j+1}\nu_{j+1}},\frac{(1+\theta_{j+1})\nu_{j+1}}{1+\theta_{j+1}\nu_{j+1}}\right),\quad \left(\frac{1+\theta_{j+1}}{1+\theta_{j+1}\nu_{j+1}},1\right).$$
(iv) The sequence $H_{q_k}(y)$, $k\in K$, converges to a number $H(y)$ where 
$$H(y)=\phi_y(g(y))=g(y).$$
\end{proposition}

\noindent{\bf Proof.}
Take $y$ from some finite interval $(\gamma_j,\gamma_{j+1})$. By definition, $q_{k+j}<q_ky<q_{k+j+1}$ for $k\in K$ big enough. It implies that $L(q_k y)=B_{k+j}$. \\
(i) We get immediately that $\mu(L(q_k y))$ tends to zero. 

\vspace{3mm}

\noindent (ii) Using equation \eqref{eq:qeta}, we get
\begin{align*}
g(y)&=\lim_{k\in K}q_k\mu(B_{k+j})=\lim_{k\in K}q_k(\eta_{k+j}+\eta_{k+j+1})=\lim_{k\in K}q_k\eta_k\left(\frac{\eta_{k+j}}{\eta_k}+\frac{\eta_{k+j+1}}{\eta_k}\right)\\
&=\frac{\delta_j+\delta_{j+1}}{\gamma_1+\delta_1}\leq \frac{1+1}{1}\leq 2.
\end{align*}
Hence, the limit $g(y)$ exists and is finite.

\vspace{3mm}

\noindent (iii) The uniform convergence of distribution functions $F_{L(q_k y)}$, $k\in\NN$, is a direct application of a result due to Coelho and de Faria (\cite{CF}, Theorem I). If we translate their result into our settings and notation, we get that if there exist limits $\nu_{j+1}$ and $\theta_{j+1}$, the distribution functions $F_{B_{k+j}}$, $k\in K$, uniformly converge to the function $\phi_y$ described in the statement of the proposition. Since $F_{B_{k+j}}=F_{L(q_ky)}$, eventually for $k\in K$, we need only to verify that the limits $\nu_{j+1}$ and $\theta_{j+1}$ exist. Since $\gamma_\ell$ exists and is finite for every $\ell\leq j+1$, the limits $\nu_\ell$ must exist and be positive, for every $\ell\leq j+1$, i.e. 
$$
\nu_{\ell}=\lim_{k\in K}[c_{k+\ell},c_{k+\ell-1},\ldots, c_{1}]>0.
$$
In particular, for every $\ell\leq j+1$, 
$$
\eta_{k+\ell-1}/\eta_{k+\ell}<c_{k+\ell+1}+1,\qquad k\in K, \ell\leq j.
$$
As
$$
\frac{\eta_{k+j}}{\eta_k}=\prod_{\ell=0}^{j-1}\frac{\eta_{k+j+1}}{\eta_{k+\ell}}
\ge\prod^{j-1}_{\ell=0}\frac1{c_{k+\ell+1}+1}>0
$$
is a lower bound for the sequence $(\eta_{k+j}/\eta_k)_{k\in K}$ we conclude that $\delta_j$ 
is strictly positive and therefore $\theta_{j+1}$ exists and is equal to $\delta_{j+1}/\delta_j$.

\vspace{3mm}

\noindent (iv) The equality $H(y)=\phi_y(g(y))$ is a direct consequence of the previous parts~(ii)
and~(iii). Since
\begin{align*}
\frac{(1+\frac{\eta_{k+j+1}}{\eta_{k+j}})\frac{q_{k+j}}{q_{k+j+1}}}{1+\frac{\eta_{k+j+1}}{\eta_{k+j}}\cdot\frac{q_{k+j}}{q_{k+j+1}}}&=\frac{(\eta_{k+j}+\eta_{k+j+1})q_{k+j}}{\eta_{k+j}q_{k+j+1}+\eta_{k+j+1} q_{k+j}}=(\eta_{k+j}+\eta_{k+j+1})q_{k+j}\\
&=\frac{(\eta_{k+j}+\eta_{k+j+1})q_{k+j}}{\eta_{k}q_{k+1}+\eta_{k+1} q_{k}}=\frac{\left(\frac{\eta_{k+j}}{\eta_k}+\frac{\eta_{k+j+1}}{\eta_k}\right)\frac{q_{k+j}}{q_k}}{\frac{q_{k+1}}{q_k}+\frac{\eta_{k+1}}{\eta_{k}}},
\end{align*}
passing to the limit on both sides yields
$$
g(y)=\frac{\delta_j+\delta_{j+1}}{\gamma_1+\delta_1}\leq \frac{(\delta_j+\delta_{j+1})\gamma_j}{\gamma_1+\delta_1}=\frac{(1+\theta_{j+1})\nu_{j+1}}{1+\theta_{j+1}\nu_{j+1}}.
$$
By the definition of $\phi_y$, we conclude that $\phi_y(g(y))=g(y)$.
\qed

\vspace{3mm}

\noindent In comparison with Theorem \ref{thm:rotations}, the last proposition does not answer the question what happens for $y\in (-\infty,1)$ and for $y\in (\gamma_j,\gamma_{j+1})$ when $\gamma_{j+1}$ is infinite. To extend the last proposition and use the formula $H(y)=\phi_y(g(y))$ also for these cases is quite complicated, because the limiting function $\phi_y$ or the limiting value $g(y)$ need not exist for every $y$ from these intervals.

%%%%%%%%%%%%%%%%%%%%%%%%%%%%%%%%%%%%%%%%
%%%%%%%%%%%  ACKNOWLEDGMENTS
%%%%%%%%%%%%%%%%%%%%%%%%%%%%%%%%%%%%%%%%
\section*{Acknowledgments}

The work on the article was initiated during the stays of authors at the University of Toulon and CPT Luminy in 2011. The authors would like to express their gratitude to the university and ANR Perturbations for the financial support of the stays. M. K. acknowledges ANR Pertubations for the financial support of the stay at CPT Luminy in 2013 and thanks to S. Vaienti for useful discussions on the topic. 

%%%%%%%%%%%%%%%%%%%%%%%%%%%%%%%%%%%%%%%%
%%%%%%%%%%%  REFERENCES
%%%%%%%%%%%%%%%%%%%%%%%%%%%%%%%%%%%%%%%%


\begin{thebibliography}{99}

\bibitem{Aba1} M Abadi: Hitting, returning and the short correlation function;
Bull.\ Braz.\ Math.\ Soc., New Series {\bf 37(4)} (2006), 1--17.

\bibitem{Col} P Collet: Statistics of closest return for some non-uniformly hyperbolic systems;
 Ergod.\ Th.\ \& Dynam.\ Syst.\ {\bf 21} (2001), 401--420

\bibitem{CF} Z Coelho and E de Faria: Limit laws of entrance times for homeomorphisms of the circle; 
Israel Jour.\ Math.\ {\bf 93} (1996), 93--112.

\bibitem{FF }A C M  Freitas and J M Freitas: On the link between dependence and independence in 
Extreme Value Theory for Dynamical Systems; Stat.\ Probab.\ Lett.\ {\bf 78} (2008) 1088--1093.

\bibitem{FF2} A C M Freitas and J M  Freitas: Extreme values for Benedicks Carleson maps; Ergod.\
 Th.\ \& Dynam.\ Sys.\ {\bf 28(4)} (2008), 1117--1133.
 
 \bibitem{FFT1} A C M Freitas, J M Freitas and M Todd: Hitting time statistics and extreme value theory;
  Probab Th.\ \& Rel.\ Fields {\bf 147(3)} (2010), 675--710.

\bibitem{FFT2} A C M Freitas, J M Freitas and M Todd: Extreme value laws in dynamical
systems for non-smooth observations;  J. Stat.\ Phys.\ {\bf 142(1)} (2011), 108--126.

\bibitem{GHN} C Gupta, M Holland and M Nicol: Extreme value theory and return time statistics 
for dispersing billiard maps and flows, Lozi maps and Lorenz-like maps;  
Ergod.\ Th.\ \& Dynam.\ Sys.\ {\bf 31(5)} (2011),  1363--1390.

\bibitem{Hay} N Haydn: Entry and return times distribution; to appear in
{\it Dynamical Systems: An International Journal dedicated to the Statistical 
Properties of Dynamical Systems}, available at {\tt http://arxiv.org/abs/1306.4475}

\bibitem{HNT} M Holland, M Nicol and A T\"or\"ok: Extreme value distributions for non-uniformly
hyperbolic dynamical systems; Trans.\ Amer.\ Math.\ Soc.\ {\bf 364(2)} (2012), 661--688. 

 \bibitem{KP} D H Kim and K K Park: The first return time properties of an irrational rotation;
  Proc.\ Amer.\ Math.\ Soc.\ {\bf 136(11)} (2008), 3941--3951. 

\bibitem{KL} M Kupsa and Y Lacroix: Asymptotics for hitting times; Annals of Probability {\bf 33(2)} (2005), 610--619.

\end{thebibliography}
\end{document}